\documentclass[12pt,reqno]{amsart}
\usepackage{amsmath,amssymb,latexsym, amsthm, amscd, mathrsfs, stmaryrd} %, txfonts}
\usepackage[linktocpage=true]{hyperref}
\hypersetup{colorlinks,linkcolor=blue,urlcolor=cyan,citecolor=blue}

\usepackage[all]{xy}
\usepackage{hyperref}
\usepackage{color}
\newcommand{\nc}{\newcommand}
\nc{\browntext}[1]{\textcolor{brown}{#1}}
\nc{\greentext}[1]{\textcolor{green}{#1}}
\nc{\redtext}[1]{\textcolor{red}{#1}}
\nc{\bluetext}[1]{\textcolor{blue}{#1}}
\nc{\brown}[1]{\browntext{ #1}}
\nc{\green}[1]{\greentext{ #1}}
\nc{\red}[1]{\redtext{ #1}}
\nc{\blue}[1]{\bluetext{ #1}}
\nc{\zb}[1]{\redtext{From zb: #1}}

\newtheorem{theorem}{Theorem}  [section]

\newtheorem{lemma}[theorem]{Lemma}

\newtheorem{proposition}[theorem]{Proposition}
\newtheorem{definition}[theorem]{Definition}

\theoremstyle{remark}

\numberwithin{equation}{section}

\newcommand{\cc}{{\mathcal C}}

\def \ch{{\mathcal H}}

\def \ct{{\mathcal T}}

\numberwithin{equation}{section}

\renewcommand{\ker}{\operatorname{Ker}\nolimits}

\newcommand{\Hom}{\operatorname{Hom}\nolimits}

\newcommand{\Aut}{\operatorname{Aut}\nolimits}
\newcommand{\Id}{\operatorname{Id}\nolimits}

\newcommand{\Ext}{\operatorname{Ext}\nolimits}

\newcommand{\mbf}{\mathbf}

\newcommand{\mrm}{\mathrm}

\newcommand{\ev}{\bar{0}}
\newcommand{\End}{\mrm{End}}

\newcommand{\N}{\mathbb N}

\newcommand{\odd}{\bar{1}}

\newcommand{\ov}{\overline}
\newcommand{\qbinom}[2]{\begin{bmatrix} #1\\#2 \end{bmatrix} }
\newcommand{\Q}{\mathbb Q}

\newcommand{\U}{\mbf U}

\newcommand{\F}{\mathbb F}

\newcommand{\arxiv}[1]{\href{http://arxiv.org/abs/#1}{\tt arXiv:\nolinkurl{#1}}}

\newcommand{\Ui}{{\mbf U}^\imath}

\newcommand{\vs}{\varsigma}

\newcommand{\Z}{\mathbb Z}

\def \fg{\mathfrak{g}}
\def \bU{{\mathbf U}}
\newcommand{\tK}{\widetilde{K}}
\def \I{\mathbb{I}}
\def \bv{v}

\newcommand{\tUi}{\widetilde{{\mathbf U}}^\imath}
\newcommand{\sqq}{{\bf v}}

\newcommand{\coker}{\operatorname{Coker}\nolimits}

\newcommand{\tU}{\widetilde{\mathbf U}}

\def \btau{{{\tau}}}
\newcommand{\tk}{\Bbbk}

\def \ff{B}

\def \cc{\mathcal C}
\def\ca{\mathcal A}
\newcommand{\Iso}{\operatorname{Iso}\nolimits}
\renewcommand{\Im}{\operatorname{Im}\nolimits}
\newcommand{\res}{\operatorname{res}\nolimits}
\newcommand{\iH}{{}^\imath\widetilde{\ch}}

\newcommand{\rep}{\operatorname{rep}\nolimits}
\newcommand{\Ker}{\operatorname{Ker}\nolimits}

\def \cI{\mathcal I}

\def \cs{{\mathcal{S}}}
\def \ch{\mathcal H}
\def \cd{\mathcal D}
\def\bfk{\mathbf{k}}
\def \ck{\mathcal K}

\def \bvs{{\boldsymbol{\varsigma}}}

\def \ce{\mathcal E}

\allowdisplaybreaks

\begin{document}

%%%%%
\title[$\imath$quantum groups of split type via derived Hall algebras]
{$\imath$quantum groups of split type via derived Hall algebras}

\author[Jiayi Chen]{Jiayi Chen}
\address{School of Mathematical Sciences,
Xiamen University, Xiamen 361005, P.R.China}
\email{jiayichen.xmu@foxmail.com}

\author[Ming Lu]{Ming Lu}
\address{Department of Mathematics, Sichuan University, Chengdu 610064, P.R.China}
\email{luming@scu.edu.cn}

\author[Shiquan Ruan]{Shiquan Ruan}
\address{ School of Mathematical Sciences,
Xiamen University, Xiamen 361005, P.R.China}
\email{sqruan@xmu.edu.cn}

%\author[Weiqiang Wang]{Weiqiang Wang}
%\address{Department of Mathematics, University of Virginia, Charlottesville, VA 22904, USA}
%\email{ww9c@virginia.edu}

\dedicatory{Dedicated to Professor Jie Xiao on the occasion of his sixtieth birthday}



\begin{abstract}
	A quantum symmetric pair consists of a quantum group $\mathbf U$ and its coideal subalgebra ${\mathbf U}^{\imath}_{\bvs}$ (called an $\imath$quantum group) with parameters $\bvs$. In this note, we use the derived Hall algebras of 1-periodic complexes to realize the $\imath$quantum groups ${\mathbf U}^{\imath}_{\bvs}$ of split type.
\end{abstract}

\maketitle
% \setcounter{tocdepth}{1}
% \tableofcontents

%\newpage
%%%%%%%
%%%%%%%

%\newpage
%%%%%%%
%%%%%%%
\section{Introduction}

\subsection{Backgrounds}

\subsubsection{Quantum groups via Hall algebras}
Inspired by Gabriel Theorem, Ringel \cite{Rin90} used the Hall algebra of representations of a Dynkin quiver to realize the positive part $\U^+$ of the quantum group, which is generalized by Green \cite{Gr95} to Kac-Moody setting. Later, Lusztig gave a geometric counterpart of Ringel's construction by considering quiver varieties, and  constructed the canonical basis of  $\U^+$.

Since then, many experts tried to realize the whole quantum groups by using Hall algebras.  There are  three typical constructions during this procedure. The first one is the realization of Kac-Moody algebras via Hall type Lie algebras of root categories given by Peng-Xiao \cite{PX00}; the second one is the Drinfeld double of Hall algebras given by Xiao \cite{X97}; the third one is the derived Hall algebras introduced by To\"{e}n \cite{T06} and Xiao-Xu \cite{XX08}.

%In 2013 \cite{Br13}, Bridgeland used the Hall algebra of 2-periodic complexes to realize
This problem was solved by Bridgeland in 2013 \cite{Br13}, who actually gave a categorical realization of the
 %realized  the
 Drinfeld double quantum groups $\tU$, a variant of $\U$ with the Cartan subalgebra doubled (with generators $K_i, K_i'$, for $i \in \I$). A  reduced version, which is the quotient of $\widetilde \U$ by the ideal generated by the central elements $K_i K_i'-1$, is then identified with $\U$.

Bridgeland's Hall algebras has found further generalizations and improvements which allow more flexibilities. Gorsky \cite{Gor18} constructed {\em semi-derived Hall algebras} for Frobenius categories. More recently, motivated by the works of Bridgeland and Gorsky, Lu-Peng \cite{LP16} formulated the {\em semi-derived (Ringel-)Hall algebras} by using $2$-periodic complexes of hereditary abelian categories. The semi-derived Ringel-Hall algebras were further constructed for arbitrary $1$-Gorenstein algebras by Lu in \cite[Appendix A]{LW19a}. %There is another geometric approach toward Bridgeland's Hall algebra developed by Qin \cite{Qin16}; cf. Scherotzke-Sibilla \cite{SS16}.

\subsubsection{$\imath$Quantum groups via $\imath$Hall algebras}

As a quantization of symmetric pairs $(\fg, \fg^\theta)$, the quantum symmetric pairs $(\U, \Ui)$ were formulated by Letzter \cite{Let99, Let02} (also cf. \cite{Ko14}) with Satake diagrams as inputs. The symmetric pairs are in bijection with the real forms of complex simple Lie algebras, according to Cartan. By definition, $\Ui =\Ui_\bvs$ is a coideal subalgebra of $\U$ depending on parameters $\bvs =(\vs_i)_{i\in \I}$ (subject to some compatibility conditions) and will be referred to as an $\imath$quantum group in this note. As suggested in \cite{BW18a}, most of the fundamental constructions in the theory of quantum groups should admit generalizations in the setting of $\imath$quantum groups; see \cite{BW18a, BK19, BW18b} for generalizations of (quasi) R-matrix and canonical bases, and also see \cite{BKLW} (and \cite{Li19}) for a geometric realization and
\cite{BSWW} for KLR type categorification of a class of (modified) $\Ui$.

Following terminologies in real group theory, we call an $\imath$quantum group {\em quasi-split} if the underlying Satake diagram does not contain any black node. In other words, the involution $\theta$ on $\fg$ is given by $\theta =\omega \circ \btau$, where $\omega$ is the Chevalley involution and $\btau$ is a diagram involution which is allowed to be $\Id$. In case $\btau=\Id$, $\Ui$ is called {\em split}. For example, a quantum group is a quasi-split $\imath$quantum group associated to the symmetric pair of diagonal type, and thus it is instructive to view $\imath$quantum groups as generalizations of quantum groups which may not admit a triangular decomposition.

In \cite{LW19a,LW20a}, Lu-Wang developed a Hall algebra approach to study the $\imath$quantum groups. More explicitly, they introduced a new kind of $1$-Gorenstein algebras: $\imath$quiver algebras, and used their $\imath$Hall algebras (aka twisted semi-derived Hall algebras of $\imath$quiver algebras) to realize the $\imath$quantum groups. In fact, analogous to Bridgeland's result, the $\imath$Hall algebra construction produces a universal $\imath$quantum group $\tUi$. The main difference between the $\imath$quantum groups $\Ui$  \`a la Letzter \cite{Let99} and the universal $\imath$quantum groups $\tUi$ (a coideal subalgebra of $\tU$) in \cite{LW19a} is that $\Ui$ depends on various parameters while $\tUi$ admits various central elements. A central reduction of $\tUi$ recovers $\Ui$.

\subsubsection{Derived Hall algebras}

Derived Hall algebras were firstly defined by To\"{e}n \cite{T06} for  DG-enhanced
triangulated categories, and then extended to arbitrary triangulated categories satisfying some homological finiteness conditions by Xiao-Xu \cite{XX08}. Unfortunately, so far it seems that derived Hall algebras could not be used to realize quantum groups, since  none of these methods can be applied to the periodic triangulated category, especially 2-periodic
triangulated category or root category, for it does not satisfy the homological finiteness conditions.
However, Xu-Chen \cite{XC13} constructed derived Hall algebras for odd periodic triangulated categories later.

A triangulated category $\ct$ is called {\em algebraic} if there exists a Frobenius category $\cc$ such that $\ct$ is equivalent to the stable category $\underline{\cc}$ of $\cc$.
For algebraic odd periodic triangulated category, Sheng-Chen-Xu \cite{SCX18} considered the relation between Gorsky's semi-derived Hall algebra of $\cc$ and Xu-Chen's derived Hall algebra of $\underline{\cc}$; compared with Gorsky's result \cite{Gor18} for any Frobenius category whose stable category is homologically finite. Lin-Peng \cite{LinP} revisited and consolidated their results by considering Lu-Peng's semi-derived Ringel-Hall algebras. Roughly speaking, derived Hall algebra of odd periodic triangulated category is a subalgebra of semi-derived Ringel-Hall algebras, and they are isomorphic by twisting certain elements (related to Cartan elements) to both kinds of Hall algebras respectively.

\subsection{Main results}

\subsubsection{Goal}
The  $\imath$Hall algebras constructed in \cite{LW19a,LW20a} are used to realize the universal $\imath$quantum groups $\tUi$ other than $\imath$quantum groups $\Ui$.
The goal of this note is to find a {\em direct} Hall algebra approach to realize Letzter's $\imath$quantum groups $\Ui$. On the other hand, people also would like to know whether derived Hall algebras could be used to realize any interesting and important quantum algebras.

We shall use derived Hall algebras of $1$-periodic derived category of quiver algebras to realize  {\em split} $\imath$quantum groups $\Ui$.

\subsubsection{Main results}

The first main result is to compare  derived Hall algebras of $1$-periodic derived categories  and Lu-Peng's semi-derived Hall algebras.

Let $\ca$ be a hereditary abelian category linearly over a finite field $\F_q$. Let $\cc_1(\ca)$ be the category of $1$-periodic complexes over $\ca$, and $\cd_1(\ca)$ the $1$-periodic derived category.
We prove in Theorem \ref{thm:main1} that the  derived Hall algebra $\cd\ch_1(\ca)$  of $\cd_1(\ca)$ is a quotient algebra of Lu-Peng's semi-derived Hall algebra  $\iH(\ca)$ by the ideal generated by the central elements $[K_\alpha]-1$ for $\alpha\in K_0(\ca)$; see \S\ref{subsec:i-Hall} for the definition of $[K_\alpha]$.
As a corollary, we give a Hall multiplication formula for $\cd\ch_1(\ca)$ by using the Hall number of $\ca$; see Corollary \ref{cor:Hall multiplication}.

Let  $\rep_\bfk^{\rm nil}( Q)$ be the category of  finite-dimensional nilpotent representations of arbitrary quiver $Q$.  Let $\cd\ch_1(\bfk Q)$ be the derived Hall algebra of $\cd_1(\rep_\bfk^{\rm nil}( Q))$.
The second main result is to realize the split $\imath$quantum group $\Ui$ by using $\cd\ch_1(\bfk Q)$; see Theorem \ref{them:main2}. In fact, Theorem \ref{them:main2} follows from Theorem \ref{thm:main1} and the realization of universal $\imath$quantum group via $\imath$Hall algebras obtained in \cite{LW19a,LW20a}.

%\bigskip
The paper is organized as follows.
Section \ref{sec:Hall} is devoted to reviewing the materials on Hall algebras, $1$-periodic complexes and semi-derived Hall algebras.
In Section \ref{sec:DHvsSDH}, we compare derived Hall algebras and semi-derived Hall algebras for the categories of $1$-periodic complexes, and Theorem \ref{thm:main1} is proved.
In Section
\ref{sec:iQGviaDH}, we review the $\imath$quantum groups and their categorical realization via semi-derived Hall algebras, and use $\cd\ch_1(\bfk Q)$ to realize the split $\imath$quantum group $\Ui$. %%, and Theorem \ref{them:main2} is proved.

\vspace{2mm}
\noindent {\bf Acknowledgements.}
We would like to thank the referee for helpful comments.
ML deeply thanks Weiqiang Wang for guiding him to study the $\imath$quantum groups, and also his continuing encouragement.
ML thanks University of Virginia for hospitality and support.  

ML is partially supported by the National Natural Science Foundation of China (No. 12171333).
SR is partially supported by the National Natural Science Foundation of China (No. 11801473) and the Fundamental Research Funds for Central Universities of China (No. 20720210006).

\section{$\imath$Hall algebras}
\label{sec:Hall}

In this paper, we take the field $\bfk=\mathbb F_q$, a finite field of $q$ elements. Let $\ca$ be a hereditary abelian category over $\bfk$. We define the $\imath$Hall algebra $\iH(\ca)$ as a twisted semi-derived Ringel-Hall algebra for the category $\cc_1(\ca)$ of $1$-periodic complexes over $\ca$.

\subsection{Hall algebras}

Let $\ce$ be an essentially small exact category in the sense of Quillen, linearly over  $\bfk=\F_q$.
%For the basics on exact categories, we refer to \cite{Buh} and references therein.
Assume that $\ce$ has finite morphism and extension spaces, i.e.,
\[
|\Hom(M,N)|<\infty,\quad |\Ext^1(M,N)|<\infty,\,\,\forall M,N\in\ce.
\]

Given objects $M,N,L\in\ce$, define $\Ext^1(M,N)_L\subseteq \Ext^1(M,N)$ as the subset parameterizing extensions whose middle term is isomorphic to $L$. We define the {\em Ringel-Hall algebra} $\ch(\ce)$ (or {\em Hall algebra} for short) to be the $\Q$-vector space whose basis is formed by the isoclasses $[M]$ of objects $M$ in $\ce$, with the multiplication defined by (see \cite{Br13})
\begin{align}
\label{eq:mult}
[M]\diamond [N]=\sum_{[L]\in \Iso(\ce)}\frac{|\Ext^1(M,N)_L|}{|\Hom(M,N)|}[L].
\end{align}
We remark that the Ringel-Hall algebra used here is the dual of the original
one defined in \cite{Rin90}.

For any three objects $L,M,N$, let
\begin{align*}
% \label{eq:Fxyz}
F_{MN}^L:= \big |\{X\subseteq L \mid X \cong N,  L/X\cong M\} \big |.
\end{align*}
The Riedtmann-Peng formula states that
\[
F_{MN}^L= \frac{|\Ext^1(M,N)_L|}{|\Hom(M,N)|} \cdot \frac{|\Aut(L)|}{|\Aut(M)| |\Aut(N)|},
\]
where $\Aut(M)$ denotes the automorphism group of $M$. For any object $M$, %let $\ell(M)$ denote the number of indecomposable direct summand of $M$, and define
let
\begin{align*}
	\label{eq:doublebrackets}
[\![M]\!]:=\frac{[M]}{|\Aut(M)|}.
\end{align*}
Then the Hall multiplication \eqref{eq:mult} can be reformulated to be
\begin{align*}
[\![M]\!]\diamond [\![N]\!]=\sum_{[\![L]\!]}F_{M,N}^L[\![L]\!],
\end{align*}
which is the version of Hall multiplication used in \cite{Rin90}.

\subsection{The category of $1$-periodic complexes}
%%%

\label{subsec:periodic}

Let $\ca$ be a hereditary abelian category  which is essentially small with finite-dimensional homomorphism and extension spaces.

 A $1$-periodic complex $X^\bullet$ in $\ca$ is a pair $(X,d)$ with $X\in\ca$ and a differential $d:X\rightarrow X$ with $d^2=0$. A morphism $(X,d) \rightarrow (Y,e)$ is given by a morphism $f:X\rightarrow Y$ in $\ca$ satisfying $f\circ d=e\circ f$. Let $\cc_1(\ca)$ be the category of all $1$-periodic complexes in $\ca$. Then $\cc_1(\ca)$ is an abelian category.
A $1$-periodic complex $X^\bullet=(X,d)$ is called acyclic if $\ker d=\Im d$. We denote by $\cc_{1,ac}(\ca)$ the full subcategory of $\cc_1(\ca)$ consisting of acyclic complexes.
Denote by $H(X^\bullet)\in\ca$ the cohomology group of $X^\bullet$, i.e., $H(X^\bullet)=\ker d/\Im d$, where $d$ is the differential of $X^\bullet$.

The category $\cc_1(\ca)$ is Frobenius with respect to the degreewise split exact structure.
The $1$-periodic homotopy category $\ck_1(\ca)$ is obtained as the stabilization
of $\cc_1(\ca)$, and the $1$-periodic derived category $\cd_1(\ca)$ is the localization of the homotopy category $\ck_1(\ca)$ with respect to quasi-isomorphisms. Both $\ck_1(\ca)$ and $\cd_1(\ca)$ are triangulated categories.

Let $\cc^b(\ca)$ be the category of bounded complexes over $\ca$ and $\cd^b(\ca)$ be the corresponding derived category with the suspension functor $\Sigma$.  Then there is a covering functor $\pi:\cc^b(\ca)\rightarrow \cc_1(\ca)$, inducing a covering functor $\pi:\cd^b(\ca)\rightarrow \cd_1(\ca)$ which is dense (see, e.g.,
\cite[Lemma 5.1]{St17}). The orbit category $\cd^b(\ca)/\Sigma$ is a triangulated category \cite{Ke2}, and we have
\begin{align*}
%\label{dereq}
\cd_1(\ca)\simeq \cd^b(\ca)/\Sigma.
\end{align*}

For any $X\in\ca$, denote the stalk complex by
\[
C_X =(X,0)
\]
 (or just by $X$ when there is no confusion), and denote by $K_X$ the following acyclic complex:
\[
K_X:=(X\oplus X, d),
\qquad \text{ where }
d=\left(\begin{array}{cc} 0&\Id \\ 0&0\end{array}\right).
\]

\begin{lemma}
[\text{\cite[Lemma 2.2]{LRW20a}}]
\label{lem:pd acyclic}
For any acyclic complex $K^\bullet$ and $p \ge 2$, we have
\begin{align*}
\label{Ext2vanish}
\Ext^p_{\cc_1(\ca)}(K^\bullet,-)=0=\Ext^p_{\cc_1(\ca)}(-,K^\bullet).
\end{align*}
\end{lemma}

For any $K^\bullet\in\cc_{1,ac}(\ca)$ and $M^\bullet\in\cc_{1}(\ca)$, by \cite[Corollary 2.4]{LRW20a}, define
\begin{align*}
\langle K^\bullet, M^\bullet\rangle=\dim_{\bfk} \Hom_{\cc_{1}(\ca)}( K^\bullet, M^\bullet)-\dim_{\bfk}\Ext^1_{\cc_{1}(\ca)}(K^\bullet, M^\bullet),
\\
\langle  M^\bullet,K^\bullet\rangle=\dim_{\bfk} \Hom_{\cc_{1}(\ca)}( M^\bullet, K^\bullet)-\dim_{\bfk}\Ext^1_{\cc_{1}(\ca)}( M^\bullet,K^\bullet).
\end{align*}
These formulas give rise to well-defined bilinear forms (called {\em Euler forms}), again denoted by $\langle \cdot, \cdot \rangle$, on the Grothendieck groups $K_0(\cc_{1,ac}(\ca))$
and $K_0(\cc_{1}(\ca))$.

Denote by $\langle \cdot,\cdot\rangle_\ca$ (also denoted by $\langle \cdot,\cdot\rangle$ if no confusion arises) the Euler form of $\ca$,  i.e.,
\[
\langle M,N\rangle_\ca=\dim_\bfk\Hom_\ca(M,N)-\dim_\bfk\Ext_\ca^1(M,N).
\]
Let $\res: \cc_1(\ca)\rightarrow\ca$ be the restriction functor by forgetting differentials.

\begin{lemma}
[\text{\cite[Lemma 2.7]{LRW20a}}]
\label{lemma compatible of Euler form}
We have
\begin{itemize}
\item[(1)]
$\langle K_X, M^\bullet\rangle = \langle X,\res (M^\bullet) \rangle_\ca$,\; $\langle M^\bullet,K_X\rangle =\langle \res(M^\bullet), X \rangle_\ca$, for $X\in \ca$, $M^\bullet\in\cc_1(\ca)$;
\item[(2)] $\langle M^\bullet,N^\bullet\rangle=\frac{1}{2}\langle \res(M^\bullet),\res(N^\bullet)\rangle_\ca$, for  $M^\bullet,N^\bullet\in\cc_{1,ac}(\ca)$.
\end{itemize}
\end{lemma}

%%%%%
\subsection{$\imath$Hall algebras}
\label{subsec:i-Hall}
%%%%%%%

Define
\[
\sqq :=\sqrt{q}.
\]

We continue to work with a hereditary abelian category $\ca$ as in \S\ref{subsec:periodic}.
Let $\ch(\cc_1(\ca))$ be the Ringel-Hall algebra of $\cc_1(\ca)$ over $\Q(\sqq)$, i.e., $\ch(\cc_1(\ca))=\bigoplus_{[M^\bullet]\in \Iso(\cc_1(\ca))} \Q(\sqq)[M^\bullet]$, with multiplication defined by
\begin{align*}
[M^\bullet]\diamond[N^\bullet]=\sum_{[L^\bullet]\in\Iso(\cc_1(\ca)) } \frac{|\Ext^1(M^\bullet,N^\bullet)_{L^\bullet}|}{|\Hom(M^\bullet,N^\bullet)|}[L^\bullet].
\end{align*}

%Let $I$ be the ideal of $\ch(\cc_1(\ca))$ generated by
%all differences $[L^\bullet]-[K^\bullet\oplus M^\bullet]$ if there is a short exact sequence $K^\bullet\rightarrowtail L^\bullet\twoheadrightarrow M^\bullet$ in $\cc_{1}(\ca)$ with $K$ acyclic.
Following \cite{LP16,LW19a,LW20a}, we consider the ideal $\cI$ of $\ch(\cc_1(\ca))$ generated by
\begin{align*}
  \label{eq:ideal}
&\Big\{ [M^\bullet]-[N^\bullet]\mid H(M^\bullet)\cong H(N^\bullet), \quad \widehat{\Im d_{M^\bullet}}=\widehat{\Im d_{N^\bullet}} \Big\}.
\end{align*}
Here we use $\widehat{X}$ to denote the image of $X\in\mathcal{A}$ in the Grothendieck group $K_0(\ca)$.
%$\widehat{\Im d_{M^\bullet}}$ and $\widehat{\Im d_{N^\bullet}}$ are respectively the images  of $\Im d_{M^\bullet}$ and $\Im d_{N^\bullet}$ in the Grothendieck group $K_0(\ca)$.
We denote
\[
\cs:=\{ a[K^\bullet] \in \ch(\cc_1(\ca))/\cI \mid a\in \Q(\sqq)^\times, K^\bullet\in \cc_1(\ca) \text{ acyclic}\},
\]
a multiplicatively closed subset of $\ch(\cc_1(\ca))/ \cI$ with the identity $[0]$.

\begin{lemma}
[\text{\cite[Proposition A.5]{LW19a}}]
The multiplicatively closed subset $\cs$ is a right Ore, right reversible subset of $\ch(\cc_1(\ca))/\cI$. Equivalently, there exists the right localization of
$\ch(\cc_1(\ca))/\cI$ with respect to $\cs$, denoted by $(\ch(\cc_1(\ca))/\cI)[\cs^{-1}]$.
\end{lemma}

The algebra $(\ch(\cc_1(\ca))/\cI)[\cs^{-1}]$ is the {\em semi-derived Ringel-Hall algebra} of $\cc_1(\ca)$ in the sense of \cite{LP16, LW19a} (also cf. \cite{Gor18}), and will be denoted by $\cs\cd\ch(\cc_1(\ca))$.

\begin{definition}[\text{\cite[Definition 2.5]{LR21}}]
  \label{def:iH}
The {\em $\imath$Hall algebra} of a hereditary abelian category $\ca$, denoted by $\iH(\ca)$,  is defined to be the twisted semi-derived Ringel-Hall algebra of $\cc_1(\ca)$, that is, the $\Q(\sqq)$-algebra on the same vector space as $\cs\cd\ch(\cc_1(\ca)) =(\ch(\cc_1(\ca))/\cI)[\cs^{-1}]$ equipped with the following modified multiplication (twisted via the restriction functor $\res: \cc_1(\ca)\rightarrow\ca$)
\begin{align}
   \label{eq:tH}
[M^\bullet]* [N^\bullet] =\sqq^{\langle \res(M^\bullet),\res(N^\bullet)\rangle_\ca} [M^\bullet]\diamond[N^\bullet].
\end{align}
\end{definition}
For any complex $M^\bullet$ and acyclic complex $K^\bullet$, we have
\[
[K^\bullet]*[M^\bullet]=[K^\bullet\oplus M^\bullet]=[M^\bullet]*[ K^\bullet].
\]

For any $\alpha\in K_0(\ca)$,  there exist $X,Y\in\ca$ such that $\alpha=\widehat{X}-\widehat{Y}$. Define $[K_\alpha]:=[K_X]* [K_Y]^{-1}$. This is well defined, see e.g.,
\cite[\S 3.2]{LP16}.
It follows that $[K_\alpha]\; (\alpha\in K_0(\ca))$ are central in the algebra $\iH(\ca)$.

%It follows that $[K_\alpha]$ are central in the algebra $\cs\cd\widetilde{\ch}(\ca)$ for any $\alpha\in K_0(\ca)$.

The {\em quantum torus} $\widetilde{\ct}(\ca)$ is defined to be the subalgebra of $\iH(\ca)$ generated by $[K_\alpha]$ for $\alpha\in K_0(\ca)$. %Then $\widetilde{\ct}(\ca)$ is the group algebra of $K_0(\ca)$.

\begin{proposition}
[\text{\cite[Proposition 2.9]{LRW20a}}]
\label{prop:hallbasis}
The folllowing hold in $\iH(\ca)$:
\begin{enumerate}
\item
The quantum torus $\widetilde{\ct}(\ca)$ is a central subalgebra of $\iH(\ca)$.
\item
The algebra $\widetilde{\ct}(\ca)$ is isomorphic to the group algebra of the abelian group $K_0(\ca)$.
\item
$\iH(\ca)$ has an ($\imath$Hall) basis given by
\begin{align*}
\{[M]*[K_\alpha]\mid [M]\in\Iso(\ca), \alpha\in K_0(\ca)\}.
\end{align*}
\end{enumerate}
\end{proposition}

%%%%%%%%%%%%%%%%%%

\section{Derived Hall algebras vs semi-derived Hall algebras}
\label{sec:DHvsSDH}

\subsection{1-Periodic derived Hall algebra}
Assume $\mathcal{A}$ is a hereditary abelian category over $\bfk$. Let $\cc_1(\ca)$ be the category of 1-periodic complexes on $\mathcal{A}$, and $\mathcal{D}_1(\mathcal{A})$ be its derived category. Observe that the isoclasses $\Iso(\mathcal{D}_1(\mathcal{A}))$ of $\mathcal{D}_1(\mathcal{A})$ coincides with the isoclasses $\Iso(\mathcal{A})$ of $\mathcal{A}$.
In the following we will identify $\Iso(\mathcal{D}_1(\mathcal{A}))$ with $\Iso(\mathcal{A})$.
\begin{lemma}\label{hom and aut in D1(A)}
For any objects $A,B\in\mathcal{A}$, we have

(1) $\Ext^1_{\mathcal{C}_1(\mathcal{A})}(A,B)\cong\Hom_{\mathcal{D}_1(\mathcal{A})}(A,B)\cong\Hom_\mathcal{A}(A,B)\oplus\Ext^1_\mathcal{A}(A,B);$

(2) $|\Aut_{\mathcal{D}_1(\mathcal{A})}(A)|= |\Aut_\mathcal{A}(A)|\cdot |\Ext^1_\mathcal{A}(A,A)|$.
\end{lemma}

\begin{proof}
(1) The formula $\Ext^1_{\mathcal{C}_1(\mathcal{A})}(A,B)\cong\Hom_{\mathcal{D}_1(\mathcal{A})}(A,B)$ is well known; see e.g. \cite{LinP} for a proof. The formula $\Hom_{\mathcal{D}_1(\mathcal{A})}(A,B)\cong\Hom_\mathcal{A}(A,B)\oplus\Ext^1_\mathcal{A}(A,B)$ follows from the fact that $\cd_1(\ca)\simeq \cd^b(\ca)/\Sigma$, where $\Sigma$ is the suspension functor of $\cd^b(\ca)$.

(2) By (1), we have $\End_{\cd_1(\ca)}(A)=\End_\ca(A)\oplus \Ext_\ca^1(A,A)$. Since $\ca$ is hereditary, it is obvious that any $f\in \Ext_\ca^1(A,A)$ is nilpotent. Thus $\Aut_{\cd_1(\ca)}(A)=\Aut_\ca(A)\oplus \Ext_\ca^1(A,A)$. Then the desired formula follows.
\end{proof}

For convenience, we denote by $$|\Aut_{\mathcal{D}_1(\mathcal{A})}(A)|=\tilde{a}_{A}\quad\text{and}\quad |\Aut_\mathcal{A}(A)|=a_{A}$$
for any $A\in \mathcal{A}$. Hence, we have $\tilde{a}_A=a_A\cdot |\Ext^1_\mathcal{A}(A,A)|$.

Following \cite{XC13}, we use the following notation for any objects $A,B,M\in\mathcal{A}$:
\[
(A,B)_M:=\{f\in\Hom_{\mathcal{D}_1(\mathcal{A})}(A,B)\mid{\rm Cone}(f)\cong M\}.
\]

\begin{lemma}\label{cardinality of triangles to extensions}
For any objects $A,B,M\in\mathcal{A}$, we have\[
|(A,B)_M|=\sum_{[X^\bullet]\in\Iso(\mathcal{C}_1(\mathcal{A})); {H(X^\bullet)}\cong M}|\Ext^1_{\mathcal{C}_1(\mathcal{A})}(A,B)_{X^\bullet}|.
\]
\end{lemma}
\begin{proof}
Observe that each exact sequence
\[
0\rightarrow B\rightarrow X^\bullet\rightarrow A\rightarrow 0
\]in $\mathcal{C}_1(\mathcal{A})$ lifts to a triangle\[
A\rightarrow B\rightarrow X^\bullet\rightarrow A
\] in $\mathcal{D}_1(\mathcal{A})$, since we have $\Sigma A\cong A$ in $\mathcal{D}_1(\mathcal{A})$.

For any two complexes $X^\bullet,Y^\bullet\in \mathcal{C}_1(\mathcal{A})$, by using $\mathcal{D}_1(\mathcal{A})=\mathcal{D}^b(\mathcal{A})/\Sigma$, we have \[
X^\bullet\cong Y^\bullet \ {\rm in}\  \mathcal{D}_1(\mathcal{A})\Longleftrightarrow H(X^\bullet)\cong H(Y^\bullet)\ {\rm in}\  \mathcal{A}.
\]
In particular, $X^\bullet\cong H(X^\bullet)$ in $\mathcal{D}_1(\mathcal{A}).$ Then the result follows from Lemma \ref{hom and aut in D1(A)} (1) immediately.
\end{proof}

Following \cite{XC13}, we denote by
\[
\{A,B\}:=\frac{1}{|\Hom_{\mathcal{D}_1(\mathcal{A})}(A,B)|}
\]
for any objects $A,B\in\mathcal{A}$. Then

\begin{lemma}\label{from {} to Ext}
For any objects $A,B\in \mathcal{A}$, we have
\begin{equation}\label{{A,B}}
\sqrt{\{A,B\}}=\sqq^{-\langle A,B\rangle_\mathcal{A}}\cdot\frac{1}{|\Ext^1_\mathcal{A}(A,B)|}.
\end{equation}
In particular,
\begin{equation}\label{{A,A}}
\sqrt{\{A,A\}}\cdot\tilde{a}_A=\sqq^{-\langle A,A\rangle_\mathcal{A}}\cdot a_A.
\end{equation}
\end{lemma}

\begin{proof}
By Lemma \ref{hom and aut in D1(A)}, we have
\begin{align*}
  \{A,B\} &= \frac{1}{|\Hom_{\mathcal{D}_1(\mathcal{A})}(A,B)|} \\
 &= \frac{1}{|\Hom_\mathcal{A}(A,B)|\cdot|\Ext^1_\mathcal{A}(A,B)|} \\
 &=q^{-\langle A,B\rangle_\mathcal{A}}\cdot \frac{1}{|\Ext^1_\mathcal{A}(A,B)|^2}.
\end{align*}
This proves \eqref{{A,B}}.
Then \eqref{{A,A}} follows from
$\tilde{a}_A=a_A\cdot |\Ext^1_\mathcal{A}(A,A)|$ immediately.
\end{proof}

For any objects $A,B,M\in\mathcal{A}$, by \cite[Corollary 2.7]{XC13} we have
\[
\frac{|(B,M)_A|}{\tilde{a}_B}\sqrt{\frac{\{B,M\}}{\{B,B\}}}=\frac{|(M,A)_B|}{\tilde{a}_A}\sqrt{\frac{\{M,A\}}{\{A,A\}}}.
\]
We denote this number by $G_{AB}^{M}$ in this paper, which satisfies the following
derived Riedtmann-Peng formula by \cite[Proposition 3.3]{SCX18}:

\begin{equation}\label{derive RP formula}
G_{AB}^M=\frac{\tilde{a}_M\cdot|(A,B)_M|}{\tilde{a}_A \tilde{a}_B}\sqrt{\frac{\{A,B\}\{M,M\}}{\{A,A\}\{B,B\}}}.
\end{equation}

Recall that $\rm{Iso}(\mathcal{D}_1(\mathcal{A}))=\rm{Iso}(\mathcal{A})$.
The \emph{1-periodic derived Hall algebra} $\mathcal{DH}_1(\mathcal{A})$ is a $\mathbb{Q}(\sqq)$-vector space with the basis $\{u_{[A]}\,|\, [A]\in {\rm{Iso}(\mathcal{D}_1(\mathcal{A}))}\}$, endowed with the multiplication defined by\[
u_{[A]}*u_{[B]}=\sum_{[M]\in \Iso(\mathcal{D}_1(\mathcal{A}))}G_{AB}^M\cdot u_{[M]}.
\]

\subsection{A homomorphism from $\imath$Hall algebra to 1-periodic derived Hall algebra}

In this subsection, we prove our first main result of this note.

Let $\ca$ be a hereditary abelian $\bfk$-linear category. Recall that $\iH(\ca)$ has an $\imath$Hall basis given by
$$\{[M]*[K_\alpha]\mid [M]\in\Iso(\ca), \alpha\in K_0(\ca)\},$$
and $[K_\alpha]\; (\alpha\in K_0(\ca))$ are central in $\iH(\ca)$.

\begin{theorem}
	\label{thm:main1}
	Let $\ca$ be a hereditary abelian $\bfk$-linear category. Then there exists an algebra epimorphism
\[
\Phi: \iH(\ca)\longrightarrow \mathcal{DH}_1(\mathcal{A}),
\]
\[
\qquad\quad\qquad [M]*[K_\alpha]\longmapsto \sqrt{\{M,M\}}\cdot \tilde{a}_M\cdot u_{[M]}
\]
with ${\rm Ker}\ \Phi=\langle [K_\alpha]-1,\alpha\in K_0(\mathcal{A})\rangle$.
\end{theorem}
\begin{proof}
For any two objects $A,B\in\mathcal{A}$, in $\iH(\ca)$ we have
\begin{align*}
 [A]*[B] &= \sqq^{{\langle A,B\rangle_\mathcal{A}}}\sum_{{[X^\bullet]}\in\Iso(\mathcal{C}_1(\mathcal{A}))}
 \frac{|\Ext^1_{\mathcal{C}_1(\mathcal{A})}(A,B)_{X^\bullet}|}{|\Hom_{\mathcal{C}_1(\mathcal{A})}(A,B)_{X^\bullet}|}
 \cdot[X^\bullet] \\
 &= \sqq^{\langle A,B\rangle_\mathcal{A}}\sum_{[X^\bullet]\in\Iso(\mathcal{C}_1(\mathcal{A}))}
 \frac{|\Ext^1_{\mathcal{C}_1(\mathcal{A})}(A,B)_{X^\bullet}|}
 {|\Hom_{\mathcal{C}_1(\mathcal{A})}(A,B)|}\cdot[H(X^\bullet)]*[K_{\Im d_{X^\bullet}}] \\
 &= \sqq^{\langle A,B\rangle_\mathcal{A}}\sum_{[M]\in \Iso(\mathcal{A})}\sum_{[X^\bullet]\in\Iso(\mathcal{C}_1(\mathcal{A})):
 H(X^\bullet)\cong M}\frac{|\Ext^1_{\mathcal{C}_1(\mathcal{A})}(A,B)_{X^\bullet}|}{|\Hom_{\mathcal{C}_1(\mathcal{A})}(A,B)|}\cdot[M]*[K_{\Im d_{X^\bullet}}],
\end{align*}
where $2\cdot\widehat{{\Im} d_{X^\bullet}}=\widehat{A}+\widehat{B}-{\widehat{H(X^\bullet)}}$.

By the definition of $\Phi$, we get\begin{align*}
 \Phi([A]*[B]) &= \sqq^{\langle A,B\rangle_\mathcal{A}}\sum_{[M]\in \Iso(\mathcal{A})}\ \sum_{[X^\bullet]\in\Iso(\mathcal{C}_1(\mathcal{A})): H(X^\bullet)\cong M}\frac{|\Ext^1_{\mathcal{C}_1(\mathcal{A})}(A,B)_{X^\bullet}|}{|\Hom_{\mathcal{C}_1(\mathcal{A})}(A,B)|}\cdot \Phi([M]) \\
 &=\frac{\sqq^{\langle A,B\rangle_\mathcal{A}}}{|\Hom_{\mathcal{A}}(A,B)|}\sum_{[M]\in \Iso(\mathcal{A})}|(A,B)_M|\cdot \Phi([M]),
% &=\sqq^{-\langle A,B\rangle_\mathcal{A}}\sum_{[M]\in \Iso(\mathcal{A})}\frac{|(A,B)_M|}{|\Ext^1_{\mathcal{A}}(A,B)|}\cdot \Phi([M]),
\end{align*}
where the second equality follows from Lemma \ref{cardinality of triangles to extensions} and the fact $\Hom_{\mathcal{C}_1(\mathcal{A})}(A,B)\cong \Hom_{\mathcal{A}}(A,B)$.

 On the other hand, we have
 \begin{align*}
 \Phi([A])*\Phi([B]) &= \sqrt{\{A,A\}}\cdot\tilde{a}_A\cdot u_{[A]}*  \sqrt{\{B,B\}}\cdot\tilde{a}_B\cdot u_{[B]}\\
 &= \sqrt{\{A,A\}\{B,B\}}\cdot\tilde{a}_A \tilde{a}_B\cdot  \sum_{[M]\in \Iso(\mathcal{A})}G_{AB}^M\cdot u_{[M]}\\
 &\stackrel{\eqref{derive RP formula}}{=} \sum_{[M]\in \Iso(\mathcal{A})}|(A,B)_M|\sqrt{\{A,B\}\{M,M\}}\cdot \tilde{a}_M\cdot u_{[M]} \\
 &= \sqrt{\{A,B\}}\cdot\sum_{[M]\in \Iso(\mathcal{A})}|(A,B)_M|\cdot \Phi([M]).
 % &\stackrel{\eqref{{A,B}}}{=} \sqq^{-\langle A,B\rangle}\sum_{[M]\in \Iso(\mathcal{A})}\frac{|(A,B)_M|}{|\Ext^1_{\mathcal{A}}(A,B)|}\cdot \Phi([M])\\
% &= \Phi([A]*[B]).
\end{align*}

By \eqref{{A,B}} we have
$$
\sqrt{\{A,B\}}=\frac{\sqq^{-\langle A,B\rangle_\mathcal{A}}}{|\Ext^1_\mathcal{A}(A,B)|}=\frac{\sqq^{\langle A,B\rangle_\mathcal{A}}}{|\Hom_\mathcal{A}(A,B)|}.
$$
It follows that $$\Phi([A]*[B])= \Phi([A]*[B]).$$

Recall that $\iH(\ca)$ has an $\imath$Hall basis given by
$\{[M]*[K_\alpha]\mid [M]\in\Iso(\ca), \alpha\in K_0(\ca)\}$,
and $[K_\alpha]\; (\alpha\in K_0(\ca))$ are central in $\iH(\ca)$. Hence $\Phi$ is an algebra homomorphism, since $\Phi([K_\alpha])=1$ for any $\alpha\in K_0(\ca)$.

Recall that $\mathcal{DH}_1(\mathcal{A})$ has a basis $\{u_{[A]}\,|\,[A]\in \rm{Iso}(\mathcal{D}_1(\mathcal{A}))\}$, and $\rm{Iso}(\mathcal{D}_1(\mathcal{A}))$ coincides with $\rm{Iso}(\mathcal{A})$.
It follows that $\Phi$ is surjective and ${\rm Ker}\ \Phi=\langle [K_\alpha]-1,\alpha\in K_0(\mathcal{A})\rangle$.
\end{proof}

As an application, we obtain the following multiplication formula for the 1-periodic derived Hall algebra $\mathcal{DH}_1(\mathcal{A})$.
\begin{proposition}
	\label{cor:Hall multiplication}
For any objects $A,B,M\in \mathcal{A}$, we have\[
u_{[A]}*u_{[B]}=\sum_{[M]\in\Iso(\mathcal{A})}\sum_{[L],[I],[N]\in \Iso(\mathcal{A})}\sqq^{\langle I,N\rangle_\mathcal{A}+\langle I,I\rangle_\mathcal{A}+\langle L,I\rangle_\mathcal{A}-\langle L,N\rangle_\mathcal{A}}\cdot \frac{a_L  a_I  a_N}{a_A  a_B}\cdot F_{NL}^MF^A_{IN}F^B_{LI}\cdot u_{[M]}.
\]
In particular,
\[
G_{AB}^M=\sum_{[L],[I],[N]\in \Iso(\mathcal{A})}\sqq^{\langle I,N\rangle_\mathcal{A}+\langle I,I\rangle_\mathcal{A}+\langle L,I\rangle_\mathcal{A}-\langle L,N\rangle_\mathcal{A}}\cdot \frac{a_L  a_I a_N}{a_A  a_B}\cdot F_{NL}^MF^A_{IN}F^B_{LI}.
\]
\end{proposition}

\begin{proof}
By \cite[Proposition 3.10]{LW20a}, we know that the following formula holds in $\iH(\ca)$:
\begin{align*}
    [A]*[B]&=
		\sum_{[L],[M],[N]\in \Iso(\mathcal{A})} \sqq^{\langle A,B\rangle_\mathcal{A}}  q^{\langle N,L\rangle_\mathcal{A} -\langle A,B\rangle_\mathcal{A}}\cdot\frac{|\Ext_\mathcal{A}^1(N, L)_{M}|}{|\Hom_\mathcal{A}(N,L)| }
		\\
		\notag
		&\cdot |\{s\in\Hom_\mathcal{A}(A,B)\mid \Ker s\cong N, \coker s\cong L		\}|\cdot [M]*[K_{\widehat{A}-\widehat{N}}].\\
\end{align*}
By Riedtmann-Peng formula we have
$$\frac{|\Ext_\mathcal{A}^1(N, L)_{M}|}{|\Hom_\mathcal{A}(N,L)| }=F_{NL}^M\cdot\frac{a_N  a_L}{a_M}.$$
Moreover, it is well-known that
$$|\{s\in\Hom_\mathcal{A}(A,B)\mid \Ker s\cong N, \coker s\cong L\}|=\sum_{[I]\in \Iso(\mathcal{A})}F_{IN}^AF_{LI}^B\cdot a_I.$$
It follows that
\begin{align*}
    [A]*[B]&=	\sum_{[M]\in \Iso(\mathcal{A})}\ \sum_{[L],[I],[N]\in \Iso(\mathcal{A})} \sqq^{-\langle A,B\rangle_\mathcal{A}} q^{\langle N,L\rangle_\mathcal{A}}\cdot F_{NL}^MF_{IN}^AF_{LI}^B\cdot\frac{a_N  a_L a_I}{a_M }		\cdot [M]*[K_{\widehat{A}-\widehat{N}}].
\end{align*}

By using Theorem \ref{thm:main1} and \eqref{{A,A}}, we get
$$\Phi([M]*[K_{\widehat{A}-\widehat{N}}])=\sqrt{\{M,M\}}\cdot \tilde{a}_M\cdot u_{[M]}=\sqq^{-\langle M,M\rangle_\mathcal{A}}\cdot a_M	\cdot u_{[M]}.$$
Hence \begin{align}
&\Phi([A]*[B])\notag\\\label{Phi(AB)}
=&\sum_{[M]\in \Iso(\mathcal{A})}\ \sum_{[L],[I],[N]\in \Iso(\mathcal{A})} \sqq^{-\langle A,B\rangle_\mathcal{A}-\langle M,M\rangle_\mathcal{A}} q^{\langle N,L\rangle_\mathcal{A}}\cdot F_{NL}^MF_{IN}^AF_{LI}^B\cdot a_N  a_L  a_I	\cdot u_{[M]}.
\end{align}
On the other hand, we have
\begin{align}
    \Phi([A]*[B])&=\Phi([A])*\Phi([B])\notag\\
    &\stackrel{\eqref{{A,A}}}{=} \sqrt{\{A,A\}}\cdot \tilde{a}_A\cdot u_{[A]} * \sqrt{\{B,B\}}\cdot \tilde{a}_B\cdot u_{[B]}\notag\\
    &=\sqq^{-\langle A,A\rangle_\mathcal{A}}\cdot a_A	\cdot u_{[A]}*\sqq^{-\langle B,B\rangle_\mathcal{A}}\cdot a_B	\cdot u_{[B]}\notag\\\label{PhiAPhiB}
    &=\sqq^{-\langle A,A\rangle_\mathcal{A}-\langle B,B\rangle_\mathcal{A}}\cdot a_A a_B\cdot\sum_{[M]\in \Iso(\mathcal{A})}G_{AB}^M\cdot u_{[M]}.
\end{align}
By comparing the coefficient of $u_{[M]}$ in \eqref{Phi(AB)} and \eqref{PhiAPhiB}, we have
\begin{equation}\label{GabMformula1}
G_{AB}^M=\sum_{[L],[I],[N]\in \Iso(\mathcal{A})} \sqq^{-\langle A,B\rangle_\mathcal{A}+\langle A,A\rangle_\mathcal{A}+\langle B,B\rangle_\mathcal{A}-\langle M,M\rangle_\mathcal{A}} q^{\langle N,L\rangle_\mathcal{A}}\cdot F_{NL}^MF_{IN}^AF_{LI}^B\cdot\frac{a_N  a_L a_I}{a_A  a_B }.
\end{equation}

Now consider the nonzero term on the right-hand side of \eqref{GabMformula1}. For $F_{IN}^AF_{LI}^B\neq 0$, we know that $A,B$ fit into the following exact sequences:
$$0\to N\to A\to I\to 0;\quad 0\to I\to B\to L\to 0.$$
Hence \begin{align}\label{ABAABB}
    &-\langle A,B\rangle_\mathcal{A}+\langle A,A\rangle_\mathcal{A}+\langle B,B\rangle_\mathcal{A}\\
    \notag &=-\langle N+I,I+L\rangle_\mathcal{A}+\langle N+I,N+I\rangle_\mathcal{A}+\langle I+L,I+L\rangle_\mathcal{A}\\\notag
    &=\langle I,N\rangle_\mathcal{A}+\langle I,I\rangle_\mathcal{A}+\langle L,I\rangle_\mathcal{A}+\langle N,N\rangle_\mathcal{A}+\langle L,L\rangle_\mathcal{A}-\langle N,L\rangle_\mathcal{A}.
\end{align}
For $F_{NL}^M\neq 0$, $M$ fits into the following exact sequence
$$0\to L\to M\to N\to 0.$$
Hence \begin{align}\label{MMNL}
    -\langle M,M\rangle_\mathcal{A}+2\langle N,L\rangle_\mathcal{A}&=    -\langle N+L,N+L\rangle_\mathcal{A}+2\langle N,L\rangle_\mathcal{A}\\\notag
    &=\langle N,L\rangle_\mathcal{A}-\langle N,N\rangle_\mathcal{A}-\langle L,L\rangle_\mathcal{A}-\langle L,N\rangle_\mathcal{A}.
\end{align}
To sum up \eqref{ABAABB} with \eqref{MMNL}, we have
\begin{align*}
    &-\langle A,B\rangle_\mathcal{A}+\langle A,A\rangle_\mathcal{A}+\langle B,B\rangle_\mathcal{A}-\langle M,M\rangle_\mathcal{A}+2\langle N,L\rangle_\mathcal{A}\\&=\langle I,N\rangle_\mathcal{A}+\langle I,I\rangle_\mathcal{A}+\langle L,I\rangle_\mathcal{A}-\langle L,N\rangle_\mathcal{A}.
\end{align*}
Therefore, \[
G_{AB}^M=\sum_{[L],[I],[N]\in \Iso(\mathcal{A})} \sqq^{\langle I,N\rangle_\mathcal{A}+\langle I,I\rangle_\mathcal{A}+\langle L,I\rangle_\mathcal{A}-\langle L,N\rangle_\mathcal{A} }\cdot F_{NL}^MF_{IN}^AF_{LI}^B\cdot\frac{a_N a_L  a_I}{a_A  a_B }.
\]
Then we finish the proof.
\end{proof}

%%%%%%%%%%%%%
\section{$\imath$quantum groups  of split type via derived Hall algebras}
\label{sec:iQGviaDH}

\subsection{$\imath$Quantum groups of split type}

Let $Q$ be a quiver (without loops) with vertex set $Q_0= \I$.
Let $n_{ij}$ be the number of edges connecting vertex $i$ and $j$. Let $C=(c_{ij})_{i,j \in \I}$ be the symmetric generalized Cartan matrix of the underlying graph of $Q$, defined by $c_{ij}=2\delta_{ij}-n_{ij}.$ Let $\fg$ be the corresponding Kac-Moody Lie algebra. Let $\alpha_i$ ($i\in\I $) be the simple roots of $\fg$.

Let $\bv$ be an indeterminant. Write $[A, B]=AB-BA$. For $r,m \in \N$, denote by
\[
[r]=\frac{\bv^r-\bv^{-r}}{\bv-\bv^{-1}},
\quad
[r]!=\prod_{i=1}^r [i], \quad \qbinom{m}{r} =\frac{[m][m-1]\ldots [m-r+1]}{[r]!}.
\]
Then $\tU := \tU_\bv(\fg)$ is defined to be the $\Q(\bv)$-algebra generated by $E_i,F_i, \tK_i,\tK_i'$, $i\in \I$, where $\tK_i, \tK_i'$ are invertible, subject to the following relations {for $i,j\in\I$}:
\begin{align}
	[E_i,F_j]= \delta_{ij} \frac{\tK_i-\tK_i'}{\bv-\bv^{-1}},  &\qquad [\tK_i,\tK_j]=[\tK_i,\tK_j']  =[\tK_i',\tK_j']=0,
	\label{eq:KK}
	\\
	\tK_i E_j=\bv^{c_{ij}} E_j \tK_i, & \qquad \tK_i F_j=\bv^{-c_{ij}} F_j \tK_i,
	\label{eq:EK}
	\\
	\tK_i' E_j=\bv^{-c_{ij}} E_j \tK_i', & \qquad \tK_i' F_j=\bv^{c_{ij}} F_j \tK_i',
	\label{eq:K2}
\end{align}
and the quantum Serre relations {for $i\neq j \in \I$},
\begin{align}
	& \sum_{r=0}^{1-c_{ij}} (-1)^r  E_i^{(r)} E_j  E_i^{(1-c_{ij}-r)}=0,
	\label{eq:serre1} \\
	& \sum_{r=0}^{1-c_{ij}} (-1)^r   F_i^{(r)} F_j  F_i^{(1-c_{ij}-r)}=0.
	\label{eq:serre2}
\end{align}
Here \[
F_i^{(n)} =F_i^n/[n]!, \quad E_i^{(n)} =E_i^n/[n]!, \quad \text{ for } n\ge 1, i\in \I.
\]
Note that $\tK_i \tK_i'$ are central in $\tU$ for all $i$.

Analogously as for $\tU$, the quantum group $\bU:=\bU_\bv(\fg)$ is defined to be the $\Q(v)$-algebra generated by $E_i,F_i, K_i, K_i^{-1}$, $i\in \I$, subject to the  relations modified from \eqref{eq:KK}--\eqref{eq:serre2} with $\tK_i$ and $\tK_i'$ replaced by $K_i$ and $K_i^{-1}$, respectively.

We define $\widetilde{\bU}^\imath$
to be the $\Q(v)$-subalgebra of $\tU$ generated by
\[
B_i= F_i +  E_{i} \tK_i',
\qquad \tk_i = \tK_i \tK_{i}', \quad \forall i \in \I.
\]
Let $\tU^{\imath 0}$ be the $\Q(v)$-subalgebra of $\tUi$ generated by $\tk_i$, for $i\in \I$.
The elements
$\tk_i$
are central in $\tUi$.

Let $\bvs=(\vs_i)\in  (\Q(\bv)^\times)^{\I}$.
Let $\Ui:=\Ui_{\bvs}$ be the $\Q(v)$-subalgebra {of $\bU$} generated by
\[
B_i= F_i+\vs_i E_{i}K_i^{-1},
\qquad  \forall i \in \I.
\]
It is known \cite{Let99, Ko14} that $\bU^\imath$ is a right coideal subalgebra of $\bU$, i.e., $\Delta (\Ui) \subset \Ui\otimes \U$; and $(\bU,\Ui)$ is called a quasi-split \emph{quantum symmetric pair}, as they specialize at $v=1$ to $(U(\fg), U(\fg^\omega))$, where $\omega$ is the Chevalley involution of $\fg$.

We call $\Ui$ an $\imath$quantum group and $\tUi$ a universal $\imath$quantum group of split type.
The algebras $\Ui:=\Ui_{\bvs}$, for $\bvs \in  (\Q(\bv)^\times)^{\I}$, are obtained from $\tUi$ by central reductions.

\begin{proposition}[\text{\cite[Propositon 6.2]{LW19a}}]
	\label{prop:QSP12}
	(1) The $\Q(v)$-algebra $\Ui$ is isomorphic to the quotient of $\tUi$ by the ideal generated by
	$\tk_i - \vs_i $ for all $i\in\I$.
	
	(2) The algebra $\widetilde{\bU}^\imath$ is a right coideal subalgebra of $\widetilde{\bU}$.
\end{proposition}

It is  well known that the $\Q(v)$-algebra $\Ui_{\bvs}$ (up to some field extension) is isomorphic to $\Ui_{\bvs^\diamond}$ for some distinguished parameters $\bvs^\diamond$ (cf. \cite{Let02}, \cite[Proposition~ 9.2]{Ko14}; also see \cite[Proposition 5]{CLW18}). Throughout this paper, we always assume the distinguished parameters to be
\begin{align*}
	\vs^\diamond_i=-v^{-2},\qquad\forall i\in\I,
\end{align*}
and $\Ui=\Ui_{\bvs^\diamond}$; see \cite[Eq. (7.1)]{LW21a}.

For the convenience of readers, we give a presentation of $\Ui$ in the following.
For $i\in \I$, generalizing the constructions in \cite{BW18a, BeW18}, we define the {\em $\imath${}divided powers} of $B_i$ to be (see also \cite{CLW18})
\begin{eqnarray*}
	&&\ff_{i,\odd}^{(m)}=\frac{1}{[m]!}\left\{ \begin{array}{ccccc} B_i\prod_{s=1}^k (B_i^2+v^{-1}[2s-1]^2 ) &  \qquad\text{if }m=2k+1, \\
		\prod_{s=1}^k (B_i^2+v^{-1}[2s-1]^2) &\text{if }m=2k; \end{array}\right.
	%\label{eq:iDPodd}
	\\
	&&\ff_{i,\ev}^{(m)}= \frac{1}{[m]!}\left\{ \begin{array}{ccccc} B_i\prod_{s=1}^k (B_i^2+v^{-1}[2s]^2 ) & \qquad\text{if }m=2k+1,\\
		\prod_{s=1}^{k} (B_i^2+v^{-1}[2s-2]^2) &\text{if }m=2k. \end{array}\right.
	%\label{eq:iDPev}
\end{eqnarray*}

%\begin{eqnarray*}
%	&&\ff_{i,\odd}^{(m)}=\frac{1}{[m]!}\left\{ \begin{array}{ccccc} B_i\prod\limits_{s=1}^k (B_i^2+v^{-1}[2s-1]^2 ) &  \qquad\text{if }m=2k+1, \\
%		\prod\limits_{s=1}^k (B_i^2+v^{-1}[2s-1]^2) &\text{if }m=2k; \end{array}\right.
%	\label{eq:iDPodd}\\\\
%	&&\ff_{i,\ev}^{(m)}= \frac{1}{[m]!}\left\{ \begin{array}{ccccc} B_i\prod\limits_{s=1}^k (B_i^2+v^{-1}[2s]^2 ) & \qquad\text{if }m=2k+1,\\
%		\prod\limits_{s=1}^{k} (B_i^2+v^{-1}[2s-2]^2) &\text{if }m=2k. \end{array}\right.
%	\label{eq:iDPev}
%\end{eqnarray*}

The following theorem is an upgrade of (and can be derived from) \cite[Theorem~3.1]{CLW18} for $\Ui$ to the setting of a universal $\imath$quantum group $\tUi$; it generalizes \cite[Proposition~6.4]{LW19a} for $\tUi$ of ADE type.

\begin{theorem}[\text{\cite[Theorem~3.1]{CLW18}}]
	\label{thm:Serre}
	Fix $\ov{p}_i\in \Z/2\Z$ for each $i\in \I$. The $\Q(v)$-algebra $\Ui$ has a presentation with generators $B_i$ $(i\in \I)$ and the relations
	\begin{align*}
			\sum_{n=0}^{1-c_{ij}} (-1)^n  B_{i, \overline{p_i}}^{(n)}B_j B_{i,\overline{c_{ij}}+\overline{p}_i}^{(1-c_{ij}-n)} &=0,\quad   \text{ if }j\neq i.
		\label{relation6}
	\end{align*}
	(This presentation is called a {\em Serre presetation} of $\Ui$.)
\end{theorem}

\subsection{$\imath$quantum groups and derived Hall algebras}
\label{subsec:iquiver algebra}

Let $Q=(Q_0,Q_1)$ be a quiver (not necessarily acyclic), and we sometimes write $\I =Q_0$. A representation $V=(V_i,V(\alpha))_{i\in Q_0,\alpha\in Q_1}$ of $Q$ is called {\em nilpotent} if for each oriented cycle $\alpha_m\cdots\alpha_1$ at a vertex $i$, the $\bfk$-linear map $V(\alpha_m)\cdots V(\alpha_1):V_i\rightarrow V_i$ is nilpotent.
Let $\rep_\bfk(Q)$ be the category  formed by finite-dimensional nilpotent representations of $Q$. Let $S_i$ be the simple modules supported at $i\in\I$.

We denote by $\iH(\bfk Q)$ the twisted semi-derived Ringel-Hall algebra of $\cc_1(\rep_\bfk^{\rm nil} (Q))$. Then we have the following result.

\begin{lemma}
	[\text{\cite{LW20a}}]
	\label{thm:main}
	Let $Q$ be an arbitrary quiver. Then there exists a $\Q(\sqq)$-algebra embedding
	\begin{align*}
		\widetilde{\psi}_Q: \tUi_{|v= \sqq} &\longrightarrow \iH(\bfk Q),
	\end{align*}
	which sends
	\begin{align*}
		B_i \mapsto \frac{-1}{q-1}[S_{i}],
		&\qquad
		\tk_i \mapsto - q^{-1}[K_{S_i}], \qquad\text{ for }i \in \I.
		%\label{eq:split}
	\end{align*}
\end{lemma}

Now we can give the second main result of this paper.
\begin{theorem}
	\label{them:main2}
		Let $Q$ be an arbitrary quiver. Then there exists a $\Q(\sqq)$-algebra embedding
	\begin{align*}
		\Psi_Q: \Ui_{|v= \sqq} &\longrightarrow \cd\ch_1(\bfk Q),
	\end{align*}
	which sends
	\begin{align*}
		B_i \mapsto -\sqq^{-1}u_{[S_{i}]},
		 \qquad\text{ for }i \in \I.
		%\label{eq:split}
	\end{align*}
\end{theorem}

\begin{proof}
From Theorem  \ref{thm:main1} and Lemma \ref{thm:main}, we have an algebra morphism:
$${\Phi\circ\widetilde{\psi}_Q}: \tUi_{|v= \sqq}\longrightarrow\cd\ch_1(\bfk Q),$$
and its kernel is $\langle \tk_i+v^{-2}\mid i\in\I\rangle$. Note that
$\Phi\circ\widetilde{\psi}_Q(B_i)=-\sqq^{-1}u_{[S_{i}]}$ for any $i\in\I$.

Furthermore, it follows from Proposition \ref{prop:QSP12} that $\Ui=\Ui_{\vs^\diamond}$ is isomorphic to the  quotient of $\tUi$ by the ideal generated by $\tk_i+v^{-2}$ for all $i\in\I$. Then $\Phi\circ\widetilde{\psi}_Q$ induces the desired $\Q(\sqq)$-algebra embedding \begin{align*}
	\Psi_Q: \Ui_{|v= \sqq} &\longrightarrow \cd\ch_1(\bfk Q).
\end{align*}
\end{proof}

%%%%%%%
%%%%%%%

\end{document}